\documentclass[10pt,twoside,reqno]{amsart}
\usepackage{bbm}
\usepackage{verbatim}
\usepackage{graphicx,xcolor}
\usepackage{enumerate}
\usepackage[all]{xy}
\usepackage[hyperindex=true,plainpages=false,colorlinks=false,pdfpagelabels]{hyperref}

\newtheorem*{theorem*}{Theorem}
\newtheorem*{remark*}{Remark}

\theoremstyle{remark}

\begin{document}

\title{Willmore spheres in the 3-sphere revisited}

\author[S. Heller]{Sebastian Heller}

\address{Institute of Differential Geometry,
Leibniz Universit\"at Hannover,
Welfengarten 1, 30167 Hannover}

\email{seb.heller@gmail.com}

\subjclass[2010]{}

\date{\today}

\begin{abstract}
Bryant \cite{Br1} classified all Willmore spheres in $3$-space to be given by minimal surfaces in $\mathbb R^3$ with embedded planar ends. This note provides new explicit formulas for genus 0 minimal surfaces in $\mathbb R^3$ with $2k+1$ embedded planar ends for all $k\geq4.$ 
Peng and Xiao claimed these examples to exist in \cite{Peng3}, 
but in the same paper they also claimed the existence of a minimal surface with 7 embedded planar ends, which was falsified by Bryant \cite{Br}. \end{abstract}

\maketitle

\section{Surfaces}
Let $\phi\colon \hat \Sigma\to S^3$ be a  compact, smooth, conformally parametrised, and immersed  surface
such that for a suitable chosen point $p\in S^3$, with its stereographic projection
$\pi_p\colon S^3\setminus\{p\}\to\mathbb R^3$,
the composition
\[f=\pi_p\circ \phi\colon \Sigma=\hat\Sigma\setminus \phi^{-1}\{p\}\to\mathbb R^3\]
is a minimal surface in $\mathbb R^3.$
We call such minimal surfaces $f$ minimal surfaces with embedded planar ends.
It was shown by Bryant \cite{Br1} that all Willmore spheres in the 3-sphere
are of this type. Conversely, every $\phi$ as above is a Willmore sphere.
By definition, (immersed) Willmore surfaces $\phi\colon \hat\Sigma\to S^3$ are the critical points for the Willmore functional
\begin{equation}\label{Willi}\mathcal W(\phi)=\int_{\hat \Sigma}(H^2-K+1)dA,\end{equation}
where $1$ is the sectional curvature of the round 3-sphere, $H$ is the  mean curvature, $dA$ and $K$ are the area form and the curvature of the induced metric of $\phi$. The Willmore energy of a Willmore sphere is given by $4 \pi (n-1), $ with $n$ being the number of ends of $f = \pi_p \circ \phi$.
\begin{remark*}
Peng and Xiao claim the existence of a minimal surface with 7 embedded planar ends in \cite{Peng3} and remark that the existence for $n = 2k+1\geq9$ follows by a long but straight forward computation.
Though the $n= 7$ case was falsified by Bryant \cite{Br}, we show that the surfaces predicted in \cite{Peng3} do exist for $n \geq 9$ by giving a simple and explicit parametrization.\end{remark*}

Let $X$ be a Riemann surface, and $g\colon X\to \mathbb K^r$  be a smooth map, where $\mathbb K\in\{\mathbb R,\mathbb C\},$ and $r\in\mathbb N$.
We denote by 
$$dg=\partial g+\bar\partial g$$
the decomposition of the differential $dg$ of $g$ into its complex linear part $\partial g$ and its complex antilinear part $\bar\partial g.$
For a minimal surface $f\colon \Sigma\to \mathbb R^3 $
there exists a holomorphic line bundle $S\to \Sigma$ with $S^2=K_\Sigma$ 
and two holomorphic sections
$s_1,s_2\in H^0(\Sigma,S)$
such that
\begin{equation}\label{wr}\partial f=(s_1^2+s_2^2,is_1^2-is_2^2,-2i s_1s_2).\end{equation}
This is called the Weierstrass representation, and the two spinors  $(s_1,s_2)$ are the Weierstrass data of $f.$
In the case a minimal surface $f\colon\Sigma\to\mathbb R^3$ with embedded planar ends its Weierstrass data  are meromorphic 
spinors  on
 $\hat\Sigma$ with first order poles at $\phi^{-1}(p)$, see for example \cite{KS} and the references therein.
Moreover, $\partial f$ has no residues at the embedded planar ends.

\subsection{Existing examples in the literature}

Peng and Xiao \cite{Peng3} considered the  following ansatz 
\[s_1=\tfrac{z (z^k-c)}{(z^k-1)(z^k-\lambda)}\sqrt{dz}\quad \text{ and }\quad s_2=\tfrac{(z^k-a) (z^k-b)}{z(z^k-1)(z^k-\lambda)}\sqrt{dz}\]
for the Weierstrass data
of a minimal surface of the $(2k+1)$-punctured sphere
\[\Sigma=\mathbb CP^1\setminus\{z\in\mathbb C\mid z(z^k-1)(z^k-\lambda)=0\},\]
where $a,b,c,\lambda\in\mathbb C$ are pairwise distinct and satisfy the following algebraic condition:
\begin{equation}\label{endeq}0=\text{res}_p{s_1^2}=\text{res}_p{s_1s_2}=\text{res}_p{s_2^2}\end{equation}
holds at every point
$p\in \{z\in\mathbb C\mid z(z^k-1)(z^k-\lambda)=0\}$. Peng and Xiao \cite{Peng3} claim that a solution $(a, b, c, \lambda)$ always exist implying via \eqref{wr} the existence of immersed minimal surfaces with $2k+1$ ends.

Appropriate solutions of \ref{endeq} for $k\in\{4,5,6\}$ and  explicit formula
for the corresponding minimal surface $f$ can be computed easily.
For example,
for $k=4,$ 
\begin{equation*}\label{abcl}
\begin{split}
a&=10 - 4 \sqrt{7} + \sqrt{\tfrac{635}{3} - 80 \sqrt{7}}, \quad b=10 + 4 \sqrt{7} + \sqrt{\tfrac{635}{3} +80 \sqrt{7}},
\quad c=-3 - 4 \sqrt{\tfrac{3}{5}}, \quad
\lambda=-31 - 8 \sqrt{15}
\end{split}
\end{equation*}
solve Equation \eqref{endeq}.
We obtain a minimal surface
$f=\Re(F)$
with 9 embedded planar ends,
where \[F\colon\mathbb CP^1\setminus\{z\in\mathbb C\mid z(z^k-1)(z^k-\lambda)=0\}\longrightarrow\mathbb C^3\] is given by
   \[F(z)=\tfrac{1}{45 z (-1 + z^4) (31 + 8 \sqrt{15}+ z^4)}
\begin{pmatrix} 6 \sqrt{-1} z^2 (-15 - 4 \sqrt{15} + 5 z^4))\\
5 (31 + 8 \sqrt{15}) - 3 z^4 (153 + 40 \sqrt{15} + 15 z^4)\\
\sqrt{-1}  (5 (31 + 8 \sqrt{15}) - 3 z^4 (147 + 40 \sqrt{15}+ 15 z^4))
   \end{pmatrix}.
   \]

\section{Curves}
The following  description of genus 0 minimal surfaces with embedded planar ends is due to Bryant \cite{Br}, see also \cite{Br3}.
Consider a genus 0 minimal surface $f\colon \mathbb CP^1\setminus\{p_1,..,p_n\}\to\mathbb R^3$ with $n$ embedded planar ends.
As $\partial f$ has no residues at the ends $p_1,..,p_n$, the surface $f$ is the real part of a meromorphic map $F\colon\mathbb CP^1\to\mathbb C^3$ with simple poles at $p_1,..,p_n.$  Since $f$ is conformally parametrised $F$ is a null curve, i.e., with respect to the standard symmetric inner product $(.,.)$ on $\mathbb C^3$ we have  $(\partial F,\partial F)=0$.
Consider $\mathbb C^5$ with the inner product \[\langle.,.\rangle=-e_0^*\otimes e_4^* -e_4^*\otimes e_0^* +e_1^*\otimes e_1^* +e_2^*\otimes e_2^* +e_3^*\otimes e_3^*,\] the 3-quadric \[\mathcal Q^3=P\{v\in\mathbb C^5\setminus\{0\}\mid \langle v,v \rangle=0\}\]
and the conformal embedding 
\[\Psi\colon\mathbb C^3\to \mathcal Q^3\,;\quad (z_1,z_2,z_3)\mapsto [\tfrac{1}{2}(z_1^2+z_2^2+z_3^2),z_1,z_2,z_3,1].\]
 For a minimal sphere $f=\Re(F)$ with $n$ embedded planar ends,
$\Psi\circ F\colon\mathbb CP^1\to \mathcal Q^3$ is an unbranched rational curve of degree $n.$ Moreover, $\Psi\circ F$ is again a null curve, i.e.,
for every local holomorphic lift $\hat\Psi$  of $\Psi\circ F$ the condition  $\langle\hat \Psi,\hat \Psi\rangle=0=\langle\partial \Psi,\partial \Psi\rangle$ holds.
Conversely,  every (nondegenerate) unbranched null curve gives rise to a minimal surface with embedded planar ends by reversing the above construction.
 
 Let $V=\mathbb C^4$ be equipped  with the 2-form
$\Omega= dx_1\wedge dx_2+dx_3\wedge d x_4.$
Consider the 5-dimensional space $W=\{\eta\in\Lambda^2V\mid \Omega(\eta)=0\}$ equipped with the non-degenerated symmetric
inner product $\tfrac{1}{2}\Omega\wedge\Omega,$ and the corresponding 3-quadric $\mathcal Q_\Omega$ of null lines in $PW$.
Identifying $(W,\tfrac{1}{2}\Omega\wedge\Omega)\cong(\mathbb C^5,\langle.,.\rangle)$ yields
 $\mathcal Q_\Omega\cong\mathcal Q^3.$
A holomorphic curve $\psi\colon\hat\Sigma\to\mathbb CP^3$  is a contact curve if   $\Omega(\hat\psi\wedge \partial\hat\psi)=0$ holds for every local holomorphic lift $\hat\psi$ of $\psi$. Then, with respect to a local holomorphic coordinate $z,$ the map 
 $z\mapsto \hat\psi\wedge\frac{\partial \hat\psi}{\partial z}$
 gives a local lift of a well-defined holomorphic curve
 \[\psi^2\colon\hat\Sigma\to \mathcal Q_\Omega\subset PW,\] the second associated curve of $\psi$. 
 It is a null curve with respect to  $\tfrac{1}{2}\Omega\wedge\Omega.$ A curve into a projective space is called nondegenerate if it is not contained in any hyperplane.
  The Klein correspondence (see \cite{Br3}) states that every nondegenerate null curve is given by  a nondegenerate contact curve in 
 $\mathbb CP^3$.
 For $\hat\Sigma=\mathbb CP^1$ and 
 $\psi$ of degree $d$ its second associated curve $\psi^2$ is unbranched if and only if
 the total branch order of  $\psi$ is $d-3$. This is a direct consequence of the Pl\"ucker relations applied to the duality between $\psi$ and its third associated curve, see \cite{Br3}. Hence, to construct genus 0 minimal surfaces with $2k+1$ embedded planar ends,  we
 have to construct rational contact curves of degree $2k$ with total branch order $2k-3.$
\subsection{Rational contact curves  of degree $2k$ with total branch order $2k-3$}
For $k\in\mathbb N\setminus\{3\}$ consider
the map $\psi\colon\mathbb CP^1\to\mathbb CP^3$ defined via the lift
\begin{equation}\label{formula}\hat\psi(z)=\begin{pmatrix}z^3\\
\tfrac{1}{6} (-6 + 13 k - 9 k^2 + 2 k^3 + 
\frac{12 (-3 + 2 k) }{-3 + k}z^k +3 k (-3 + 2 k) z^k
   - 6 z^{2 k})\\
\tfrac{1}{2} (-3 + 2 k) z (-2 + k + 2 z^k)\\
z^2 (-1 + k + z^k)\end{pmatrix}.\end{equation}
This is a nondegenerate rational curve of degree $2k$ if $k\in\mathbb N^{>3}$. It can be directly verified that it is a contact curve, i.e.,
$\Omega(\hat\psi\wedge\tfrac{\partial \hat\psi}{\partial z})=0.$
For $k\in\mathbb N^{>3}$ its branch points are  at the $k$-th roots of unity and at $z=\infty$. The branch order at the roots of unity is 1, and at $z=\infty$ the branch order is $k-3.$
Hence,  the total branch order is $2k-3.$ The construction fails for $k\leq3:$
for $k=1$, the curve 
is of degree $2$ and  branched at $z=0$ and $z=1.$ Thus, its second associate curve cannot be unbranched.
For $k=2$, the degree of the curve is  $2\neq4.$ For $k=3$ and $k=0$, the formula \eqref{formula} gives a point in $\mathbb CP^3$. For $k=4$ we obtain the first valid example
\[\hat\psi(z)=(z^3,5+20 z^4-z^8,5z+5 z^5,3z^2+z^6).\]
Together with the  examples  of unbranched null curves of even degree $d\geq4$ and the nonexistence results of Bryant in \cite{Br},  this shows:
\begin{theorem*}
There exist a Willmore sphere $\phi_n\colon \mathbb CP^1\to S^3$ with Willmore energy $4\pi (n-1)$
if and only if  $n\in\mathbb N\setminus\{2,3,5,7\}.$
\end{theorem*}
For completeness, we state the formula for the genus 0 minimal surfaces $f=\Re(F)$ with $(2k+1)$ embedded planar ends corresponding to the contact curve \eqref{formula}.
 Note that these surfaces are only determined up to Goursat transformations \cite{Goursat}, see also \cite{Br}.  We leave it as an exercise for the interested reader to verify that we have reobtained the surfaces of
  Peng and Xiao \cite{Peng3} up to a Goursat transformation and reparametrisation. The meromorphic map $F\colon\mathbb CP^1\to \mathbb C^3$  is given by
  \begin{equation*}
  F(z)=
  \begin{pmatrix}
  \tfrac{  -12 \sqrt{-1} (-3 + k) (-3 + 2 k) z^2 (2 - 3 k + k^2 - 2 z^k)}{4 z ((-3 + k) (-1 + k)^2 (-3 + 2 k) - 6 (-1 + k) (-3 + 2 k) z^k - 
   3 (-3 + k) z^{2 k})}\\
  \tfrac{ (-1 + k) (-12 (3 - 2 k)^2 (-2 + k) z^k - 
   12 (-3 + k) (-3 + 2 k) z^{
    2 k} + (-3 + k) ((6 - 7 k + 2 k^2)^2 - 12 z^4))}{4 z ((-3 + k) (-1 + k)^2 (-3 + 2 k) - 6 (-1 + k) (-3 + 2 k) z^k - 
   3 (-3 + k) z^{2 k})}\\
   \tfrac{\sqrt{-1} (-1 + k) (12 (3 - 2 k)^2 (-2 + k) z^k + 
   12 (-3 + k) (-3 + 2 k) z^{
    2 k} + (-3 + k) (-(6 - 7 k + 2 k^2)^2 - 12 z^4))}{4 z ((-3 + k) (-1 + k)^2 (-3 + 2 k) - 6 (-1 + k) (-3 + 2 k) z^k - 
   3 (-3 + k) z^{2 k})}\\
  \end{pmatrix}.
  \end{equation*}

\end{document}